\documentclass[12pt,a4paper]{article}
\usepackage{amstext}
\usepackage{fancyhdr}
\usepackage{amsfonts,graphicx,bezier, amssymb}
\usepackage{amsmath}
\usepackage{caption}
\usepackage{mathtools}
\usepackage[utf8]{inputenc}
\usepackage{tikz}
\usetikzlibrary{arrows.meta}
\newenvironment{prf}{\noindent{\bf{Proof:}}~~}{\hfill\rule{1ex}{1ex}\vskip1.5ex}
\newcommand{\Z}{\mathbb Z}

\newcommand{\N}{\mathbb N}

\newcommand{\beqa}{\begin{eqnarray}}
\newcommand{\enqa}{\end{eqnarray}}
\newcommand{\beq}{\begin{eqnarray*}}
\newcommand{\enq}{\end{eqnarray*}}

\newtheorem{rem}{Remark}[section]

\newtheorem{defn}{Definition}[section]

\newtheorem{thm}{Theorem}[section]

\newtheorem{lem}{Lemma}[section]

\makeatletter
\providecommand*{\twoheadrightarrowfill@}{%
  \arrowfill@\relbar\relbar\twoheadrightarrow
}
\providecommand*{\twoheadleftarrowfill@}{%
  \arrowfill@\twoheadleftarrow\relbar\relbar
}
\providecommand*{\xtwoheadrightarrow}[2][]{%
  \ext@arrow 0579\twoheadrightarrowfill@{#1}{#2}%
}
\providecommand*{\xtwoheadleftarrow}[2][]{%
  \ext@arrow 5097\twoheadleftarrowfill@{#1}{#2}%
}
\makeatother

\begin{document}

\begin{center}
{\bf\Large A Cousin Complex for the Quantum Projective Space}
\end{center}

\begin{center}

		\begin{center}
Kobi Kremnizer\footnote{Mathematical Institute, University of Oxford, Andrew Wiles Building, Radcliffe
Observatory Quarter, Woodstock Road, Oxford, OX2 6GG, England}  and  	David Ssevviiri\footnote{Department of Mathematics, Makerere University,  P.O. BOX 7062, Kampala, Uganda}, \footnote{Corresponding author}

E-mail: Yakov.Kremnitzer@maths.ox.ac.uk and  david.ssevviiri@mak.ac.ug 
			
		\end{center}

\end{center}

\begin{abstract}
 Grothendieck constructed    a Cousin complex  for abelian sheaves on an arbitrary topological space.  In a special setting, its dual called the BGG resolution is applicable in representation theory. 
  Arkhipov proposed a complex whose dual is only  suitable for representation theory of quantum groups at roots of unity of prime order. It is desirable to get one which works for quantum groups at all roots of unity. For a quantum projective space, we  provide such a complex.
 \end{abstract}

{\bf Keywords}: Noncommutative geometry, Quantum groups, Cousin complex
\vspace*{0.4cm}

{\bf MSC 2010} Mathematics Subject Classification: 20G42, 81R50, 16E05, 16T20, 18G10, 51E15

 \section{Introduction}
 
 \paragraph\noindent

 To every topological space $X$  with a filtration $X=Z_0\supseteq Z_1 \supseteq Z_2 \supseteq \cdots$ of closed subsets
  together with an Abelian sheaf on $X$,  one 
 can associate a complex in terms of local cohomologies called the Grothendieck-Cousin complex, see for instance \cite[p. 235]{Ha} and 
 \cite[p. 352]{Ke}. If $A$ is a finitely generated commutative algebra over a field $k$, then there is a canonical complex $\mathcal{K}_A$
of $A$-modules, called the residue complex. $\mathcal{K}_A$ is characterized as the Cousin
complex of the twisted inverse image $\pi^{!}k$, where $\pi : X = \text{Spec} (A) \rightarrow k$ is
the structural morphism, see \cite{Ha} and \cite{Ye}.

\paragraph\noindent

 Let $R$ be a commutative Noetherian ring. By taking $X:=\text{Spec}(R)$, the spectrum of  $R$,
 Sharp in \cite{Sh} formulated an algebraic  version of this complex in terms of modules $M$ defined over  $R$.
 This complex $C_R(M)$ characterizes   
   Gorenstein rings, Cohen-Macaulay modules as well as balanced big Cohen-Macaulay modules.  
   A ring 
 $R$ is Gorenstein if and only if the complex $C_R(R)$ provides a minimal injective resolution of $R$, \cite[Theorem 5.4]{Sh}.  
 A nonzero finitely generated $R$-module $M$ is Cohen-Macaulay if and only if the Cousin complex $C_R(M)$ is exact, \cite[Theorem 2.4]{Sh3}. 
 For the characterization of balanced big Cohen-Macaulay modules, see \cite[Theorem 4.1]{Sh2}.

 \paragraph\noindent
 
 Let $\mathcal{B}$ be the variety
 of Borel subgroups of a reductive algebraic group $G$. If $X:=\mathcal{B}$ and the filtration of $X$ is given by the Schubert varieties,
 then the Cousin complex associated to this filtration is exact and dual to the so called  Bernstein-Gelfand-Gelfand  (BGG)  resolution in representation theory, see for instance; \cite[Sec 1(a)]{Falk} and \cite[Section IX]{Ku}. Yekutieli and Zhang in \cite{YZ} gave a Cousin complex
  on  a derived category of modules defined over a not necessarily commutative ring and used it to define a dualising 
  complex for non-commutative rings.

  \paragraph\noindent
  
  A question of interest is to define stratifications and their associated Cousin complexes in the non-commutative setting.
This can then be  used to  construct the dual BGG resolution (also called the contragredient BGG resolution) for quantum groups. A candidate for the dual BGG resolution for quantum
 groups at roots of unity was proposed by Arkhipov in \cite{Ar1} and \cite{Ar2}. He proved that these are resolutions in case of a root of unity of prime order
  by reducing $\text{mod}~p$ and using Kempf's dual BGG complex. It is desired to have  a construction which works for all roots of unity.  
  The aforementioned quantum dual BGG resolution  will have applications in the representation theory of quantum groups 
   at roots of unity.

   \paragraph\noindent

  In this paper, we achieve the desired goal but in some special case. We
  construct in terms of graded modules modulo torsion submodules,
  a Cousin complex for the quantum projective space $\mathbb{P}^n_q$. $\mathbb{P}^n_q$ is
  the non-commutative analog of the   projective space $\mathbb{P}^n$ defined on the commutative polynomial rings. It turns out that 
  when the deformation $q=1$, the obtained complex is the already well known Sharp-Cousin complex on modules over commutative rings obtained by 
  using a filtration on   the topological space $\mathbb{P}^n := \text{Spec}(k[X_1, \cdots, X_{n+1}])$, where $k$  is a field.
  
  \paragraph\noindent
  
  The main ideas behind our proofs are generally those used by Kempf in \cite{Ke} for the Cousin complex. The setting in which
  we use them is however different. Instead of sheaves, we have graded modules defined over quantized algebras modulo torsion modules.
  We have a quantum projective space $\mathbb{P}^n_q$ in the place of a topological space $X$.   The second key aspect about the methods used in the paper is the definition of a 
   quantum global section supported at a quantum closed subset, see Definition \ref{d}. Lastly, we utilize the noncommutative
   analogue of  Artin-Serre's Theorem as set out by Artin and Zhang in \cite{AZ}, in particular, see Equivalence (\ref{eq4}).

   \section{Main Results}
 
 \paragraph\noindent
 Let $A$ be a Noetherian graded $k$-algebra. By $\text{ \bf Gr}(A)$, we denote a category of graded $A$-modules and by 
 $\text{ \bf Tors}(A)$ a subcategory of $\text{ \bf Gr}(A)$ consisting of  torsion modules. An element $m$ of an $A$-module $M$ is {\it torsion} if there exists an integer $N$ such that $A_{\geq N}m=0$ where $A_{\geq N}:=\oplus_{n\geq N} A_n$. A module is {\it torsion} if all its elements are torsion elements. By Serre's Theorem, there  is an 
 equivalence between the category of quasi-coherent sheaves on the projective space $\text{\bf{Proj}}(A)$ and the category of 
 graded $A$-modules modulo torsion modules. We write
 
 \begin{equation}\label{eq1}
  \text{\bf{Qcoh}}(\text{\bf{Proj}}(A))\cong  \text{\bf Gr}(A)/\text{\bf Tors}(A).
 \end{equation}

 Since we are interested in the quantized space $\mathbb{P}^n_q$, we use the projective space $\mathbb{P}^n(k)$ over a field  $k$. So, we have
 
  \begin{equation}\label{eq2}
  \text{\bf{Qcoh}}(\mathbb{P}^n(k))\cong  \text{\bf Gr}(A)/\text{\bf Tors}(A).
 \end{equation}
 
 \paragraph\noindent
 
 The algebra $A$ in this special case is nothing but a polynomial ring over $k$ in $n+1$ indeterminates.
 By utilizing the view point of Artin and Zhang \cite{AZ}, we pass to the noncommutative setting. We now instead have 
 the quantized $n$-projective space $\mathbb{P}^n_q$ and the quantized  polynomial rings   $A_{q, n}$. 
 
 \begin{equation}\label{eq3} 
A_{q, n}:=k\langle x_1, x_2, \cdots x_{n+1}\rangle /\langle x_ix_j=qx_jx_i\rangle, i\not=j. 
\end{equation}

 We are now led to the equivalence

  \begin{equation}\label{eq4}
  \text{\bf{Qcoh}}(\mathbb{P}^n_q(k))\cong  \text{\bf Gr}(A_{q, n})/\text{\bf Tors}(A_{q, n}).
 \end{equation}
 
 We denote the category $\text{\bf Gr}(A_{q, n})/\text{\bf Tors}(A_{q, n})$ by $\text{\bf Tails}(A_{q, n})$.   
     The filtration $$\mathbb{P}^n\supseteq  \mathbb{P}^{n-1}\supseteq \cdots\supseteq \mathbb{P}^1\supseteq\mathbb{P}^0=\{\infty\}$$ 
     suggests a non-existent       quantum filtration 
     $$\mathbb{P}^n_q\supseteq  \mathbb{P}^{n-1}_q\supseteq \cdots\supseteq \mathbb{P}^1_q\supseteq\mathbb{P}^0_q$$ which is 
     respectively in a bijective      correspondence with

     $$\text{\bf Tails}(A_{q, n}) \supseteq \text{\bf Tails}(A_{q, n-1}) \supseteq \cdots \supseteq      \text{\bf Tails}(A_{q, 1}) \supseteq \text{\bf Tails}(A_{q, 0})$$
     
     of graded modules. This  correspondence motivates Definition \ref{defn1}.

    \subsection{The quantum global section}  
     
     \begin{defn}\label{defn1} \rm Let $k\text{-{\bf Vect}}$ denote the category of $k$-vector spaces. For integers  $0\leq j \leq n$,      
      let $\Gamma$ be the functor  $$\Gamma~ : ~\text{\bf{Qcoh}}(\mathbb{P}^j_q)\longrightarrow k\text{-{\bf Vect}}$$
      $$ M=\oplus_{\lambda\in \Lambda}M_{\lambda}\longmapsto \Gamma(M)=M_0,$$ where $M_0$ is the degree zero homogeneous component of the 
      graded  $A_{q, n}$-module $M$.
     \end{defn}
     
     \paragraph\noindent

    The functor $\Gamma$ is left exact and the category $\text{\bf{Qcoh}}(\mathbb{P}^n_q)$ is abelian and has enough injectives.  We can therefore 
     pass to the derived category by taking right derived  functors of $\Gamma$, i.e., we have 
     
     \begin{equation}\label{eq6}
      R\Gamma ~:~ \text{\bf{D}}(\text{\bf{Qcoh}}(\mathbb{P}^j_q))\longrightarrow \text{\bf{D}}(k\text{\bf-Vect})
     \end{equation}

for all integers $ 0\leq j \leq n$.

\paragraph\noindent

Since the affine $n$-space $\mathbb{A}^n = \mathbb{P}^n\setminus\mathbb{P}^{n-1}$, we let
$\mathbb{A}^n_q := \mathbb{P}^n_q\setminus\mathbb{P}^{n-1}_q$
be the non-existent quantum open subset of $\mathbb{P}^n_q$. We then have

$$ \text{\bf{Qcoh}}(\mathbb{A}^n_q) = \text{\bf{Qcoh}}( \mathbb{P}^n_q\setminus 
\mathbb{P}^{n-1}_q) $$ 
$$\cong\text{ graded}~(A_{q, n}[x_{n+1}^{-1}])\text{-modules modulo torsion modules}$$
$$=\text{graded}\left(k\langle x_1, x_2, \cdots, x_n\rangle/\langle x_ix_j -q x_jx_i\rangle\right)\text{-modules modulo torsion modules},
$$ for $i\not=j$. 

\paragraph\noindent

If $M\in \text{\bf Tails}(A_{q, n})$, we already know from Definition \ref{defn1} that $\Gamma(M)=M_0$. Now define 

\begin{equation}\label{eq7}
\Gamma(\mathbb{A}^n_q, M)= \Gamma\left(\mathbb{P}^n_q\setminus \mathbb{P}^{n-1}_q, M\right):= \left(M[x_{n+1}^{-1}]\right)_0.
\end{equation}

We define $\Gamma(\mathbb{P}_q^n, M)$ supported at $\mathbb{P}_q^{n-1}$ as 

\begin{equation}\label{eq8}
 \Gamma_{\mathbb{P}^{n-1}_q}(\mathbb{P}^n_q, M):=
 ~\text{Ker}\left(\Gamma(\mathbb{P}_q^n, M) \longrightarrow\Gamma(\mathbb{P}^n_q\setminus \mathbb{P}^{n-1}_q, M ) \right)
                                                = ~\text{Ker}\left(M_0\longrightarrow  \left(M[x_{n+1}^{-1}]\right)_0 \right).
 \end{equation}

 \paragraph\noindent
 In general,  we have
 
 \begin{defn}\label{d}\rm The quantum global section  $\Gamma(-, M)$   is given by
  
 \begin{enumerate}
  \item  $\Gamma(\mathbb{P}^n_q, M):= M_0$;
  \item  for $0<i \leq n$, $\Gamma(\mathbb{A}^{n-i+1}_q, M) := \Gamma(\mathbb{P}^n_q\setminus \mathbb{P}^{n-i}_q, M)=\left(M[x_{n+1}^{-1}, 
             x^{-1}_n, \cdots, x_{n-i+2}^{-1}]\right)_0$;
  \item  for $0<i \leq n$, $\Gamma_{\mathbb{P}^{n-i}_q}(\mathbb{P}^n_q, M):=\text{Ker}\left(\Gamma(\mathbb{P}^{n}_q, M\right)
             \longrightarrow \Gamma(\mathbb{A}^{n-i+1}_q, M))\\
             ~~~~~~~~~~~~~~~~~~~~~~~~~~~~~~~~~~~~~~=\text{Ker}\left(M_0\longrightarrow  \left(M[x_{n+1}^{-1}, 
             x^{-1}_n, \cdots, x_{n-i+2}^{-1}]\right)_0 \right).$
 \end{enumerate}

 \end{defn}

 \begin{lem}\label{l1}
  For any module $M\in \text{\bf Tails}(A_{q, n})$, we have
  \begin{enumerate}
   \item $\Gamma_{\mathbb{P}^n_q}(\mathbb{P}^n_q, M)= \Gamma(\mathbb{P}^n_q, M)$, 
   \item $\Gamma_{\emptyset}(\mathbb{P}^n_q, M)= 0$,
   \item for all $0< i \leq n$,  
   $$0\longrightarrow \Gamma_{\mathbb{P}^{n-i}_q}(\mathbb{P}^n_q, M)\longrightarrow \Gamma(\mathbb{P}^n_q, M)\longrightarrow
     \Gamma(\mathbb{A}^{n-i+1}_q, M)\longrightarrow 0$$ is a short exact sequence.
  \end{enumerate}
 \end{lem}

 \begin{prf}
  The lemma is immediate from the definition of $\Gamma(\mathbb{P}^n_q, M)$ supported at $\mathbb{P}^n_q$, $\emptyset$ and $\mathbb{P}_q^{n-i}$
  respectively.
 \end{prf}

 \begin{defn}\rm\label{dd}
  
 For  a filtration $\mathbb{P}^0\subseteq \mathbb{P}^1\subseteq   \cdots \subseteq \mathbb{P}^n $ and $0\leq i \leq n$, we define 
 
 \begin{equation}  
\Gamma_{(\mathbb{P}^{n-i}_q) / ({\mathbb{P}_q^{n-(i+1)}})} (\mathbb{P}^n_q, M):=
  \Gamma_{\mathbb{P}_q^{n-i}}(\mathbb{P}^n_q ,M)/\Gamma_{\mathbb{P}^{n-(i+1)}_q}(\mathbb{P}^n_q, M)
 \end{equation}
 which is the quotient module 
           \begin{equation}
             \frac{\text{Ker}\left(M_0\longrightarrow  \left(M[x_{n+1}^{-1}, 
             x^{-1}_n, \cdots, x_{n-i+2}^{-1}]\right)_0 \right) }          
             {\text{Ker}\left(M_0\longrightarrow  \left(M[x_{n+1}^{-1}, 
             x^{-1}_n, \cdots, x_{n-i+1}^{-1}]\right)_0 \right)} .
           \end{equation}
 \end{defn}

 \begin{lem}\label{l2} Let  $M\in \text{\bf Tails}(A_{q, n})$.
 
 \begin{enumerate}
  \item For any integers $0<z_2\leq z_1 \leq n$,  $$0\rightarrow \Gamma_{\mathbb{P}^{z_2}_q}(\mathbb{P}^n_q, M)\rightarrow
  \Gamma_{\mathbb{P}^{z_1}_q}(\mathbb{P}^n_q, M)
  \rightarrow\Gamma_{{(\mathbb{P}^{z_1}_q})/({\mathbb{P}^{z_2}_q})}(\mathbb{P}^n_q, M)\rightarrow 0$$ is a short exact sequence.
  \item $ \Gamma_{(\mathbb{P}^{n}_q) /({\emptyset})}(\mathbb{P}^n_q, M)=
        \Gamma_{\mathbb{P}^{n}_q}(\mathbb{P}^n_q, M)$.
  \item  $\Gamma_{({\mathbb{P}^{n}_q)/({\mathbb{P}^{n}_q})}}(\mathbb{P}^n_q, M)=0$.
  \item  for any integers   $0 < z_2\leq z_1 \leq n$ and $0 < w_2\leq w_1 \leq n$, such that $z_1 \leq w_1$ and $z_2 \leq w_2$, 
  there is a $\Z$-graded module homomorphism 
   
  $$\Gamma_{({\mathbb{P}^{z_1}_q)/({\mathbb{P}^{z_2}_q})}}(\mathbb{P}^n_q, M)\rightarrow 
   \Gamma_{({\mathbb{P}^{w_1}_q)/({\mathbb{P}^{w_2}_q})}}(\mathbb{P}^n_q, M).$$
   
   \item  for any integers   $0 <z_3\leq z_2\leq z_1 \leq n$    
  $$0 \rightarrow \Gamma_{({\mathbb{P}^{z_2}_q)/({\mathbb{P}^{z_3}_q})}}(\mathbb{P}^n_q, M)\rightarrow 
   \Gamma_{({\mathbb{P}^{z_1}_q)/({\mathbb{P}^{z_3}_q})}}(\mathbb{P}^n_q, M)
   \rightarrow     \Gamma_{({\mathbb{P}^{z_1}_q)/({\mathbb{P}^{z_2}_q})}}(\mathbb{P}^n_q, M) \rightarrow 0$$
    is a short exact sequence.
 \end{enumerate}
   \end{lem}
   
   \begin{prf}
    1), 2) and 3) are immediate from Definition \ref{dd} and the definition of $\Gamma(\mathbb{P}^n_q, M)$ at a support. 4) is due to the fact that quotients of graded
     modules are graded.       5) is a consequence of 4) and the  isomorphism theorem.
   \end{prf}

 \begin{lem}\label{l3} If  $M\in \text{\bf Tails}(A_{q, n})$ and $0< z_2\leq z_1\leq n$ are integers, then there is 
 a natural injection 
       
        \begin{equation}   \label{ne}         
       \Gamma_{(\mathbb{P}_q^{z_1})/(\mathbb{P}_q^{z_2})}(\mathbb{P}^n_{q}, M) \hookrightarrow 
  \Gamma_{(\mathbb{P}_q^{z_1})\setminus(\mathbb{P}_q^{z_2})}(\mathbb{P}^n_{q}\setminus \mathbb{P}_q^{z_2}, M) 
        \end{equation}

   which becomes an isomorphism whenever the restriction 
   $$\gamma: \Gamma (\mathbb{P}^n_{q}, M)  \rightarrow   \Gamma(\mathbb{P}^n_{q}\setminus \mathbb{P}_q^{z_2}, M)$$ is surjective.
 \end{lem}

 \begin{prf}
  Remember that Inclusion (\ref{ne}) can be re-written as  
  $$ \frac{\text{Ker}\left(\Gamma(\mathbb{P}_q^n, M)\rightarrow \Gamma(\mathbb{P}_q^n\setminus \mathbb{P}_q^{z_1}, M)\right)}
  {\text{Ker}\left(\Gamma(\mathbb{P}_q^n, M)\rightarrow \Gamma(\mathbb{P}_q^n\setminus \mathbb{P}_q^{z_2}, M)\right)} \hookrightarrow
  \text{Ker}\left(\Gamma(\mathbb{P}_q^n\setminus \mathbb{P}_q^{z_2}, M)\rightarrow \Gamma(\mathbb{P}_q^{z_1}\setminus
   \mathbb{P}_q^{z_2}, M)\right). $$   
   The restriction $\gamma$ takes $\Gamma_{\mathbb{P}_q^{z_1}}(\mathbb{P}_q^{n}, M)$ into 
   $\Gamma_{\mathbb{P}_q^{z_1}\setminus \mathbb{P}_q^{z_2}}(\mathbb{P}_q^n\setminus \mathbb{P}_q^{z_2}, M)$. 
   The kernel of  $\gamma$ is $\Gamma_{\mathbb{P}_q^{z_2}}(\mathbb{P}_q^n, M)$. $\gamma$ therefore induces an injection of 
    $\frac{\Gamma_{\mathbb{P}_q^{z_1}}(\mathbb{P}_q^n, M)}{\Gamma_{\mathbb{P}_q^{z_2}}(\mathbb{P}_q^n, M)}=
    \Gamma_{(\mathbb{P}_q^{z_1})/(\mathbb{P}_q^{z_2})}(\mathbb{P}^n_{q}, M)$ into 
    $\Gamma_{(\mathbb{P}_q^{z_1})\setminus(\mathbb{P}_q^{z_2})}(\mathbb{P}^n_{q}\setminus \mathbb{P}_q^{z_2}, M).$
    For the second part, take any map $f$  in $\Gamma(\mathbb{P}_q^n\setminus \mathbb{P}_q^{z_2}, M)$  which has support in 
    $(\mathbb{P}_q^{z_1})\setminus(\mathbb{P}_q^{z_2})$. By the assumption, $f$ can be extended to $f'$ in 
    $\Gamma(\mathbb{P}_q^n, M)$. As $f$ and $f'$ have the same restriction to  $(\mathbb{P}_q^{z_1})\setminus(\mathbb{P}_q^{z_2})$, 
    $f'$ must have support in $\mathbb{P}_q^{z_1}$. It follows that in this case, $ \Gamma_{\mathbb{P}_q^{z_1}}(\mathbb{P}^n_{q}, M)$
    maps surjectively onto $\Gamma_{(\mathbb{P}_q^{z_1})\setminus(\mathbb{P}_q^{z_2})}(\mathbb{P}^n_{q}\setminus \mathbb{P}_q^{z_2}, M) $
     which completes the proof.
 \end{prf}

 \paragraph\noindent

 By following \cite[Page 250]{AZ},  \cite[Page 413]{BK} and \cite[Definition 7.1]{Kee}, we define an ample autoequivalence.

 \begin{defn}\rm \label{ample}  Let  $A_{q, n}$  be the quantized $k$-algebra of $A_n$,  and $s$ an autoequivalence on the abelian category $\text{\bf Tails}(A_{n, q})$. 
 $s$ is said to be {\it ample} if:
  \begin{itemize}
   \item[(B1)] for all $M\in \text{\bf Tails}(A_{n, q})$, there exists  positive integers $l_1, \cdots , l_p$ and an epimorphism 
             $$ \bigoplus_{i=1}^p s^{-l_i}(A_{n, q})\twoheadrightarrow M;$$
   \item[(B2)] for all epimorphisms $M\twoheadrightarrow N$, with $M, N\in \text{ \bf Tails}(A_{n, q})$ there exists $t_0$ such that for all $t\geq t_0$ and $1\leq k \leq n$
   $$ \Gamma_{\mathbb{P}_q^{k}}(\mathbb{P}_q^{n}, s^t(M))\rightarrow \Gamma_{\mathbb{P}_q^{k}}(\mathbb{P}_q^{n}, s^t(N))$$ is an epimorphism. 
  \end{itemize}

 \end{defn}
 
 \begin{lem}\label{autoequivalence} Let $M\in \text{\bf Tails}(A_{q, n})$. 
 The operation $$s^t(M):=\left(M[x_{n+1}^{-1}, 
             x^{-1}_n, \cdots, x_{n-t+2}^{-1}]\right)_0 $$
defines an autoequivalence  on the abelian category $\text{\bf Tails}(A_{q, n})$.
 \end{lem}

 \begin{lem}\label{l4} Let  $M_1, M_2, M_3$ and $M$ be modules in $\text{\bf Tails}(A_{q, n})$. If  $s$ is  an ample autoequivalence on $\text{\bf {Tails}}(A_{n, q})$, then there exists $t_0\in \N$ such that for all $t\geq t_0$, the functors $\Gamma(\mathbb{P}^n_{q}, s^t(-))$, $\Gamma_{\mathbb{P}_q^{z_1}}(\mathbb{P}^n_{q}, s^t(-))$ and 
           $\Gamma_{(\mathbb{P}_q^{z_1})/(\mathbb{P}_q^{z_2})}(\mathbb{P}^n_{q}, s^t(-))$ 
  preserve  the short exact sequence $$0\rightarrow M_1 \rightarrow M_2 \rightarrow M_3\rightarrow 0.$$

 
 \end{lem}

 \begin{prf}
   Since the three functors are left exact, it is enough to show that, they preserve epimorphisms for $t\geq t_0$. However, this becomes immediate by the hypothesis of the autoequivalence $s$ being ample, in particular (B2) of Definition \ref{autoequivalence}.    
 \end{prf}

  \subsection{Passing to the derived functor}
 \paragraph\noindent
 Since the category $\text{\bf Tails}(A_{q, n})$ is abelian and has enough injectives, for each $M\in \text{\bf Tails}(A_{q, n})$, we define a  cohomology module by  $$H^i_{\mathbb{P}^j_q}(\mathbb{P}^n_q, M):= H^i{\bf R}\Gamma_{\mathbb{P}^j_q}(\mathbb{P}^n_q, M)$$
 for $0 < j \leq n$; the cohomologies of the complex obtained by taking  the right derived functor of $\Gamma_{\mathbb{P}^j_q}(\mathbb{P}^n_q, M)$.

 \begin{lem}\label{l5}
  If $M\in \text{\bf Tails}(A_{q, n})$ and $s$ is  an ample autoequivalence on $\text{\bf Tails}(A_{q, n})$ such that $0< z_2 \leq z_1\leq n$, then there
  exists $t_0\in \N$ such that for all $t\geq t_0$,

     \begin{equation}\label{emu}     
\Gamma_{(\mathbb{P}_q^{z_1})/(\mathbb{P}_q^{z_2})}(\mathbb{P}^n_{q}, s^t(M)) \cong
  H^0_{(\mathbb{P}_q^{z_1})/(\mathbb{P}_q^{z_2})}(\mathbb{P}^n_{q}, s^t(M))   
     \end{equation}
     
  and  $~ \text{for all }~i>0$
  
  \begin{equation}\label{bbiri}
   H^i_{(\mathbb{P}_q^{z_1})/(\mathbb{P}_q^{z_2})}(\mathbb{P}^n_{q},s^t(M))=0.
  \end{equation}

   \end{lem}

 \begin{prf} Isomorphism (\ref{emu}) is a standard result about left exact functors and the associated cohomology module at degree zero.  
 The functor $\Gamma_{(\mathbb{P}_q^{z_1})/(\mathbb{P}_q^{z_2})}(\mathbb{P}^n_{q}, -)$  is left exact. However, since $s$ is an ample autoequivalence, there exists $t\geq t_0$ such that it preserves epimorphisms and hence it becomes exact. This leads to Equation (\ref{bbiri}).
 \end{prf}

 \begin{lem}\label{l6}
  Let  $0\rightarrow M_1 \rightarrow M_2\rightarrow M_3\rightarrow 0$  be a short exact sequence of    modules in
   the category   $\text{\bf Tails}(A_{q, n})$. If $0< z_2\leq z_1\leq n$ are integers and $s$ is an ample autoequivalence on $\text{\bf Tails}(A_{q, n})$, then 
   there   exists $t_0\in \N$ such that for all $t\geq t_0$,  we    have a long exact sequence   
   $$ 0 \rightarrow H^0_{(\mathbb{P}_q^{z_1})/(\mathbb{P}_q^{z_2})}(\mathbb{P}^n_{q}, s^t(M_1)) 
   \rightarrow H^0_{(\mathbb{P}_q^{z_1})/(\mathbb{P}_q^{z_2})}(\mathbb{P}^n_{q}, s^t(M_2)) 
   \rightarrow H^0_{(\mathbb{P}_q^{z_1})/(\mathbb{P}_q^{z_2})}(\mathbb{P}^n_{q}, s^t(M_3)) $$
   
   $$ \stackrel{\delta}{\rightarrow}     H^1_{(\mathbb{P}_q^{z_1})/(\mathbb{P}_q^{z_2})}(\mathbb{P}^n_{q}, s^t(M_1)) \rightarrow 
   H^1_{(\mathbb{P}_q^{z_1})/(\mathbb{P}_q^{z_2})}(\mathbb{P}^n_{q}, s^t(M_2)) \rightarrow \cdots$$
 \end{lem}

 \begin{prf}
By  Lemma \ref{l4}, $$ 0 \rightarrow \Gamma_{(\mathbb{P}_q^{z_1})/(\mathbb{P}_q^{z_2})}(\mathbb{P}^n_{q}, s^t(M_1)) 
   \rightarrow \Gamma_{(\mathbb{P}_q^{z_1})/(\mathbb{P}_q^{z_2})}(\mathbb{P}^n_{q}, s^t(M_2)) 
   \rightarrow \Gamma_{(\mathbb{P}_q^{z_1})/(\mathbb{P}_q^{z_2})}(\mathbb{P}^n_{q}, s^t(M_3))\rightarrow 0 $$ is a short exact sequence which by (\ref{emu}) is isomorphic to
   $$ 0 \rightarrow H^0_{(\mathbb{P}_q^{z_1})/(\mathbb{P}_q^{z_2})}(\mathbb{P}^n_{q}, s^t(M_1)) 
   \rightarrow H^0_{(\mathbb{P}_q^{z_1})/(\mathbb{P}_q^{z_2})}(\mathbb{P}^n_{q}, s^t(M_2)) 
   \rightarrow H^0_{(\mathbb{P}_q^{z_1})/(\mathbb{P}_q^{z_2})}(\mathbb{P}^n_{q}, s^t(M_3))\rightarrow 0 .$$ The long exact sequence obtained with connecting maps $\delta$ is a standard result in homological algebra.
\end{prf}
 
 \begin{lem}\label{l7}
 Let $M\in \text{\bf Tails}(A_{q, n})$. If 
   $z_1, z_2, z_3$  are integers such that $0 <z_3\leq z_2\leq z_1 \leq n$, then 
 
  $$0 \rightarrow H^0_{({\mathbb{P}^{z_2}_q)/({\mathbb{P}^{z_3}_q})}}(\mathbb{P}^n_q, M)\rightarrow 
   H^0_{({\mathbb{P}^{z_1}_q)/({\mathbb{P}^{z_3}_q})}}(\mathbb{P}^n_q, M)
   \rightarrow     H^0_{({\mathbb{P}^{z_1}_q)/({\mathbb{P}^{z_2}_q})}}(\mathbb{P}^n_q, M)  $$
   $$
   \stackrel{\delta}{\rightarrow} H^1_{({\mathbb{P}^{z_2}_q)/({\mathbb{P}^{z_3}_q})}}(\mathbb{P}^n_q, M)\rightarrow
   H^1_{({\mathbb{P}^{z_1}_q)/({\mathbb{P}^{z_3}_q})}}(\mathbb{P}^n_q, M) \rightarrow \cdots$$
    is a long exact sequence.
 \end{lem}

 \begin{prf}
  First apply Lemma \ref{l2}, part 5).  Just like in the proof of Lemma \ref{l6}, 
  the resulting long exact sequence with connecting maps $\delta$ is   well-known.
 \end{prf}

 \paragraph\noindent
 As in Lemma \ref{l7}, it is easy to see that for any three integers $z_{i+2}<z_{i+1}<z_i$, there exists a boundary map 
 
 $$H^i_{({\mathbb{P}^{z_i}_q)/({\mathbb{P}^{z_{i+1}}_q})}}(\mathbb{P}^n_q, M)\stackrel{\delta_i}\rightarrow
   H^{i+1}_{({\mathbb{P}^{z_{i+1}}_q)/({\mathbb{P}^{z_{i+2}}_q})}}(\mathbb{P}^n_q, M) . $$ This leads us to a complex given in Theorem \ref{thm1} where the cohomology at each point is also given in terms of isomorphic modules.

 \begin{thm}\label{thm1}
  For any $M\in \text{\bf Tails}(A_{q, n})$
  and integers  $n=z_0\geq z_1 \geq z_2 \geq \cdots$,

  $$0 \rightarrow \Gamma(\mathbb{P}^n_q, M) \stackrel{e}{\rightarrow} H^0_{ (\mathbb{P}_q^{z_0})/({\mathbb{P}_q^{z_1}})}(\mathbb{P}^n_q, M)
  \stackrel{d_0}{\rightarrow} H^1_{(\mathbb{P}_q^{z_1})/(\mathbb{P}_q^{z_2})}(\mathbb{P}^n_q, M)
  \stackrel{d_1}{\rightarrow} H^2_{(\mathbb{P}_q^{z_2})/(\mathbb{P}_q^{z_3})}(\mathbb{P}^n_q, M) \stackrel{d_2}{\rightarrow}\cdots$$
  is a complex which we   call the   Cousin complex of the quantum projective space $\mathbb{P}^n_q$.  
  Moreover, 
  
  \begin{enumerate}
   \item the kernel of $e$ is $\Gamma_{\mathbb{P}_q^{z_1}}(\mathbb{P}_q^n, M)$;
   \item the kernel of $d_0$ modulo the image of $e$ is isomorphic to the quotient
   $$\frac{\left(H^0_{(\mathbb{P}_q^n)/(\mathbb{P}_q^{z_2})}(\mathbb{P}_q^n, M)\right)}
   {\left( \text{Im}(e') +  H^0_{(\mathbb{P}_q^{z_1})/(\mathbb{P}_q^{z_2})}(\mathbb{P}_q^n, M)   \right)},$$ where
      
  $$e': \Gamma(\mathbb{P}_q^n, M) \rightarrow H^0_{(\mathbb{P}_q^n)/(\mathbb{P}_q^{z_2})}(\mathbb{P}_q^n, M)$$ is the natural homomorphism;
  
  \item if $i>0$, the kernel of $d_i$ modulo the image of $d_{i-1}$ is isomorphic to the image of $$H^i_{(\mathbb{P}_q^{z_i})/(\mathbb{P}_q^{z_{i+2}})}(\mathbb{P}_q^n, M)\rightarrow H^i_{(\mathbb{P}_q^{z_{i-1}})/(\mathbb{P}_q^{z_{i+1}})}(\mathbb{P}_q^n, M). $$
  \end{enumerate}

 \end{thm}

 \begin{prf} For brevity, just like in \cite[Lemma 7.8]{Ke}, in this proof, we drop $(\mathbb{P}_q^n, M)$ which occurs at the end of all symbols involved.
 \begin{enumerate}
  \item  Since $\mathbb{P}_q^n = \mathbb{P}_q^{z_0}$, the natural homomorphism $\Gamma\rightarrow H^0_{{\mathbb{P}_q^{z_0}}/{\mathbb{P}_q^{z_i}}}$ exists for all $i$. By the equality 
  $\Gamma_{\mathbb{P}_q^{n}}  = H^0_{{\mathbb{P}_q^{z_0}}/{\emptyset}}$ and Lemma \ref{l7}, the kernel of 
  $\Gamma\rightarrow H^0_{{\mathbb{P}_q^{z_0}}/{\mathbb{P}_q^{z_i}}}$ is $\Gamma_{\mathbb{P}_q^{z_i}}$.
  \item  Consider the commutative diagram in Figure \ref{cdiag}.
    
  \begin{figure}[ht] 
  
  \begin{center}
   
  \begin{tikzpicture}[scale=1.1]
 \node at (-9, 1) { 0};  
\draw[->, line width = 0.2mm] (-8.9, 1) -- (-8.2, 1); 
 \node at (-7.3, 1) { $H^0_{{\mathbb{P}_q^{z_1}}/{\mathbb{P}_q^{z_2}}}$}; 
 
\draw[->, line width = 0.2mm] (-6.5, 1) -- (-5.3, 1); 
 \node at (-4.3, 1) { $H^0_{{\mathbb{P}_q^{z_0}}/{\mathbb{P}_q^{z_2}}}$}; 
 
\draw[->, line width = 0.2mm] (-3.5, 1) -- (-2.2, 1); 

 \node at (-1.2, 1) { $H^0_{{\mathbb{P}_q^{z_0}}/{\mathbb{P}_q^{z_1}}}$};  
 
\draw[->, line width = 0.2mm] (-0.3, 1) -- (0.7, 1); 
 \node at (1.5, 1) { $H^1_{{\mathbb{P}_q^{z_1}}/{\mathbb{P}_q^{z_2}}}$};  

 \node at (-4.5, 3) { $\Gamma$};  
 \node at (-1.5, 3) { $\Gamma$};  
 \node at (-4.8, 2.2) { $e'$};
\draw[->, line width = 0.2mm] (-4.5, 2.6) -- (-4.5, 1.5); 
\draw[->, line width = 0.2mm] (-1.5, 2.6) -- (-1.5, 1.5); 
\draw[->, line width = 0.2mm] (-4.2, 3) -- (-1.9, 3); 
\node at (0.3, 1.3) { $d_0$}; 

\node at (-1.15, 2.1) { $e$}; 

\node at (-3.1, 3.3) { $\cong$}; 
 \end{tikzpicture}
 \end{center}
\caption{}\label{cdiag}   
  \end{figure}
  
  By Lemma \ref{l7}, the bottom row of the commutative diagram in Figure \ref{cdiag} is exact. So, $d_0\circ e=0$. This establishes 2.

  \begin{figure}[ht] 
  
  \begin{center}

\begin{tikzpicture}[scale=1.3]

 \node at (-4.5, 1) { $H^j_{{\mathbb{P}_q^{z_i}}/{\mathbb{P}_q^{z_{i+2}}}}$}; 
 
\draw[->, line width = 0.2mm] (-3.5, 1) -- (-2.2, 1); 

 \node at (-1.2, 1) { $H^j_{{\mathbb{P}_q^{z_i}}/{\mathbb{P}_q^{z_{i+1}}}}$};  
 
\draw[->, line width = 0.2mm] (-0.3, 1) -- (1.1, 1); 
 \node at (2.2, 1) { $H^{j+1}_{{\mathbb{P}_q^{z_{i+1}}}/{\mathbb{P}_q^{z_{i+2}}}}$};

 \node at (-1.2, -1.5) { $H^{j}_{{\mathbb{P}_q^{z_{i-1}}}/{\mathbb{P}_q^{z_{i+1}}}}$};

 \node at (-4.5, 3) {$H^{j-1}_{\mathbb{P}_q^{z_i-1}/{\mathbb{P}_q^{z_i}}}$};  
 \node at (-1.2, 3) { $H^{j-1}_{{\mathbb{P}_q^{z_i-1}}/{\mathbb{P}_q^{z_i}}}$};  
\draw[->, line width = 0.2mm] (-4.5, 2.6) -- (-4.5, 1.5); 
\draw[->, line width = 0.2mm] (-1.5, 2.6) -- (-1.5, 1.5); 
\draw[->, line width = 0.2mm] (-3.7, 3) -- (-2.30, 3);

\draw[->, line width = 0.2mm] (-1.5, 0.5) -- (-1.5, -1.0); 
\node at (0.3, 1.3) { $d_i$};

\node at (-1.0, 2.1) { $d_{i-1}$};

\node at (-3.1, 3.3) { $\cong$};  
 \end{tikzpicture}
 \end{center}
\caption{}\label{cd2}   
  \end{figure}

  It remains to prove that  $d_i\circ d_{i-1}=0$ for all $i>0$ and that statement (3) also holds. Let $j$ be a positive integer and consider the commutative diagram in Figure \ref{cd2}. Lemma \ref{l7} establishes the existence of this diagram and the exactness of its rows and columns. If $j=i$, then $d_i$ and $d_{i-1}$ are the homomorphisms in the Cousin complex given in  Theorem \ref{thm1}. From the commutativity and exactness of the diagram, $d_i\circ d_{i-1}=0$.
 \end{enumerate}

 \end{prf}

  \begin{rem}\rm
   
 We remark that the local cohomology modules that appear in the Cousin complex of the quantum projective space given in Theorem \ref{thm1} are actually modules for the quantum group with divided powers and are also comodules for the algebra of functions on the quantum Borel, see for instance; \cite{Hu} and \cite{Ar3} respectively.
 
  \end{rem}

 \section*{Acknowledgement}
 \paragraph\noindent
 Part of this work was  done  while the second author was visiting the University of Oxford with 
 support from the Africa-Oxford 
 initiative. He is grateful to  the first author and Prof. Balazs Szendroi for  the hospitality.


\addcontentsline{toc}{chapter}{Bibliography}
  \begin{thebibliography}{99}
  
  \bibitem{Ar1} Sergey M. Arkhipov,  Algebraic construction of contragradient
quasi-Verma modules in positive characteristic,
  {\it Representation Theory of Algebraic Groups and Quantum Groups}, 
  Advanced Studies in Pure Mathematics, {\bf 40}, (2004), 27-68. 
  
  \bibitem{Ar2} Sergey M.  Arkhipov,   Semiinfinite cohomology of contragradient Weyl modules over small quantum groups, Preprint math.QA/9906071 (1999),
1-20.




\bibitem{Ar3} Sergey Arkhipov and Dennis Gaitsgory, Another realization of the category of modules over the small quantum group, 
{\it Adv. Math.}, {\bf 173}, (2003), 114–143.

  
\bibitem{AZ} Michael Artin and James J. Zhang, noncommutative projective schemes, {\it Adv.  Math.}, {\bf 109}, (1994), 228--287.
  
  
\bibitem{BK} Erik Backelin and  Kobi Kremnizer, Quantum flag varieties, equivariant quantum D-modules, and localization of quantum groups, {\it Adv.  Math.}, {\bf 203}(2), (2006), 408-429.

\bibitem{Falk} Michael Falk, Vadim Schechtman, and Alexander Varchenko, BGG resolutions via configuration spaces, {\it J. \'{E}c. polytech. Math.},  {\bf 1} (2014), 225--245.
 
 \bibitem{Ha} Robin Hartshorne, Residues and duality, Berlin-Heidelberg, New York,  Springer, Lecture Notes in Mathematics No: 20, 1966.
 
 \bibitem{Hu} Naihong Hu, Quantum divided power algebra, $q$-derivatives and some new quantum groups, {\it J. Algebra}, {\bf 232}, (2000), 507--540.
  
 \bibitem{Kee} Dennis R. Keeler, Ample filters of invertible sheaves, {\it J. Algebra}, {\bf 259}(1), (2003), 243--283.
 
 \bibitem{Ke} George Kempf, The Grothendieck-Cousin complex of an induced representation, {\it Adv. Math.}, {\bf 29}, (1978),
                     310--396.
                     
 \bibitem{Ku} Shrawan Kumar, Kac-Moody groups, their flag varieties and representation theory, Springer Science and Business Media,  2002. 
 
\bibitem{Sh2} Rodney Y. Sharp, A Cousin complex characterization of balanced big Cohen-Macaulay modules, 
           {\it Quart.  J. Math.}, {\bf 33}(4), 1982, 471--485.
   \bibitem{Sh3} Rodney Y. Sharp, Gorenstein modules, {\it Math. Z.}, {\bf 115}, (1970), 117--139.
 
                       
   \bibitem{Sh} Rodney Y. Sharp, The Cousin complex for a module over a commutative Noetherian ring, {\it Math. Z}, {\bf 112}, (1969), 340--356.
   
\bibitem{Ye} Amnon Yekutieli, An explicit construction of the Grothendieck residue complex (with an appendix by P. Sastry) {\it  Asterisque}, {\bf 208}, (1992).
                     
  \bibitem{YZ} Amnon Yekutieli and James J. Zhang, Residue complexes over noncommutative rings, {\it J. Algebra}, {\bf 259}, (2003), 451--493.
  
 
 \end{thebibliography}
\end{document}